\def\A{\leavevmode\setbox0\hbox{A}\lower1.4ex\hbox to\wd0
{\hss`}\kern-.9\wd0A}
\def\E{\leavevmode\setbox0\hbox{E}\lower1.4ex\hbox to\wd0
{\hss`\/}\kern-.9\wd0E}
\def\a{\leavevmode\setbox0\hbox{a}\lower1.4ex\hbox to\wd0
{\hss`\/}\kern-\wd0a}
\def\e{\leavevmode\setbox0\hbox{e}\lower1.4ex\hbox to\wd0
{\hss`\/}\kern-\wd0e}
\newcommand{\be}{\begin{equation}}
\newcommand{\ee}{\end{equation}}
\newcommand{\ba}{\begin{array}}
\newcommand{\ea}{\end{array}}
\newcommand{\beqn}{\begin{eqnarray}}
\newcommand{\eeqn}{\end{eqnarray}}
\newcommand{\zero}{\setcounter{equation}{0} \par}
\font\symb=msam7
\def\znakr{\raise1.5pt\hbox{\symb\char66\kern-2pt\char74}}
\def\znakl{\raise1.5pt\hbox{\symb\char73\kern-2pt\char67}}
\def\normalsize{
\setlength{\textheight}{23cm}
\setlength{\textwidth}{15cm}
\setlength{\topmargin}{-2.0cm}
\setlength{\hoffset}{-0.5cm}
\setlength{\leftmargin}{-1cm}
\setlength{\rightmargin}{2.0cm}}
\documentstyle[11pt]{article}
\normalsize

\begin{document}
\baselineskip=20pt
\title{The Two-Dimensional Quantum\\ Galilei Groups 
\vspace{4ex}}
\author{Emil Kowalczyk \thanks{Supported by
\L{}\'od\'z University Grant no 580} \\
Department of Theoretical Physics II \\
University of \L \'od\'z \\
ul.Pomorska 149/153, 90-236 \L \'od\'z, Poland\\
e-mail: emilkow@krysia.uni.lodz.pl
\vspace{3ex}}
\date{ }
\setcounter{section}{0}
\setcounter{page}{1}
\maketitle
\begin{abstract}
The Poisson structures on two-dimensional Galilei
group, classified in the author previous paper are
quantized. The dual quantum Galilei Lie algebras are
found.
\end{abstract}
\newpage
\section{Introduction}
	In the present paper we continue the study of deformed 
nonrelativistic symmetries. In the previous paper [1] 
all Lie-Poisson structures on two-dimensional Galilei 
group were classified up to automorphism.
Below we quantize these structures showing that  the consistent 
Hopf algebras are obtained. We find also the corresponding quantum Lie 
algebras by straightforward application of duality rules.

	As a result we obtain two families of quantum groups and
quantum Lie algebras, one depends on two and the other depends
on three parameters. Various limiting cases appear after sending 
appropriate subsets of parameters to infinity.
\zero
\section{Poisson Structures on Two-Dimensional\\ Galilei Group}
	Recently all Lie bialgebra structures on two-dimensional 
Galilei algebra have been found and their Lie-Poisson counterparts 
 have been classified [1].
It appeared that, up to the automorphisms, there are nine inequivalent 
bialgebra
structures on two-dimensional Galilean Lie algebra (see table I 
in Ref.[1] ).
The corresponding Lie-Poisson structures on two-dimensional Galilei 
Group
 are shown in table II in [1].

	As  it is seen from this table, in order to impose Lie-Poisson 
structures on Galilei Group two dimensionful constants are needed. 
They can attain arbitrary nonzero values, different choice being related 
by automorphisms.
The only revelant free parameter is the dimensionless parameter
 $ \varepsilon $;
different values of $ \varepsilon $  correspond to nonequivalent Lie-Poisson
structures.

	It is worth to note that this relatively rich family of 
nonequivalent 
Lie-Poisson structures contains only one coboundary. It is in a sharp 
contrast 
with semisimple case [2] as well as the case of four-dimensional 
Poincare group [3].

The Lie-Poisson structures described provide the the starting point for 
obtaining
two-dimensional quantum Galilei groups. These groups will be here 
constructed by
applying the naive quantization procedure consisting in replacing the
Poisson brackets by commutators (and supplying the resulting 
commutation rules
with imaginary unit and appropriate dimensionful constants). It is 
obvious
from the table II in Ref.[1] that no ordering problems can appear.

	As it was mentioned, the above classification of Lie-Poisson 
structures is
 complete up to the automorphisms. However, this can not be apriori 
taken for granted
in the quantum case due to the noncommutativity of generators. 
This phenomenon
is well known in quantum mechanics:
not every cannonical transformation can be lifted to the unitary one.

We do not  attempt here to classify all nonequivalent quantum structures; 
rather
we find the quantum counterparts of "cannonical" Poisson structures 
described
in table II in [1]. 
\zero
\section{Two-dimensional quantum Galilei groups}
We apply the standard quantization procedure to the Poisson
structures given in Ref.[1]. The result can be 
summarised
as follows.There are two families of quantum groups, one depending on 
two and the other depending on three dimensionful parameters.
The relevant commutation rules read, respectively:
\[ [a,v]={-i\over2 \kappa} v^2 \hspace{7em} \]
\[ 
[a, \tau]={i\over \kappa} a + {i\over \rho} v \hspace{4em} (A)
\]
\[
[v, \tau]={i\over \kappa} v 
\hspace{8em} \]
 and
\[
[a,v]={{i\over \alpha} \tau}-{i\over \sigma}v \hspace{6em} 
\]
\[
[a, \tau]={{i\over \sigma} \tau}+{i\over \lambda}v \hspace{4em} (B)
\]
\[
[v, \tau]=0 \hspace{9em}
.\]
	The algebra (A) corresponds to the case 1 of table II in Ref.[1]
while the algebra 
(B) to all remaining
structures. (Some of them can be obtained by taking an appriopriate 
parameters to infinity).

The dimensions of constans $ \alpha ,\kappa ,\rho ,\sigma, \lambda, $
 are as follows: 
\be [ \alpha]={{s^2}\over{m^2}} \hspace{2em} 
 [ \kappa]={1\over s}   \hspace{2em}
 [ \lambda]={1\over s^2} \hspace{2em}  
 [ \rho]={1\over s^2}  \hspace{2em} 
 [ \sigma]={1\over m.} 
\ee 

The commutation rules (A) and (B) are supplied with the standard 
coproduct,
antipode, counit and *- structures,
\[ \vspace{-1ex}
\Delta (a)=a\otimes I + I\otimes a + v\otimes \tau  
\] 
\be \vspace{-1ex}
S(a)= -a + v \tau  \hspace{6em}
\ee 
\[ \vspace{-1ex}
\varepsilon (a)=0 \hspace{9em}
\]
\[ \vspace{-1ex}
a^* = a, \hspace{9em}
\]
\[
\Delta (v)=v\otimes I + I\otimes v \hspace{4em}  
\]
\be
S(v)= -v   \hspace{8em}
\ee
\[
\varepsilon (v)=0 \hspace{9em}
\]
\[
v^* = v, \hspace{9em}
\]
\[
\Delta ( \tau)= \tau\otimes I + I\otimes \tau \hspace{4em}  
\]
\be
S( \tau)= - \tau  \hspace{8em}
\ee
\[
\varepsilon ( \tau)=0 \hspace{9em}
\]
\[
 \tau^* = \tau. \hspace{9em}
\]

	We have checked that all relations providing our commutation 
rules with the 
structure of *-Hopf algebra are fulfiled.
Therefore we obtain two Hopf algebra structures ( (A) and (B) ).
\zero
\section{Duality and quantum Lie algebras}

	Having found the quantum Galilei groups one can ask what is 
the structure
of their dual Hopf algebras, i.e. the quantum Lie algebras.

In the present section we find them by straightforward aplication
 of duality 
rules.
It is well known that the dual Hopf algebra can be defined by the 
following
duality rules
\be
<XY, \Phi>=< X \otimes Y, \Delta \Phi>
\ee
\be
<X, \Phi \Psi>=< \Delta X, \Phi \otimes \Psi>;
\ee 
also the *-structure can be defined by the formulae [4]
\be
<X^*,\Phi>=<X,S^{-1}(\Phi^*)>
\ee
provided the following identity holds
\be
S^{-1} (\Phi)=[S(\Phi^*)]^* .
\ee
It is easy to check that the equation (4.4)  
is in our case fulfiled. 

In order to find explicit form of quantum Lie algebras we used the 
following scheme [4].
First, we define the Lie algebra generators by adopting the classical 
duality relations 
\be
<X, \Phi>=-i{d\over dt} \Phi(e^{itX})\mid _{t=0} ,
\ee
$\vspace{-1ex}$ i.e.$ \vspace{-1ex}$
\be 
<H, \tau^k a^l v^m>= i \delta_{1k} \delta_{0m} \delta_{0l} \hspace{2em}
\ee
\be
<P, \tau^k a^l v^m>= i \delta_{0k} \delta_{0m} \delta_{1l} \hspace{2em}
\ee
\be
<K, \tau^k a^l v^m>= i \delta_{0k} \delta_{1m} \delta_{0l}. \hspace{2em}
\ee

	These rules can be compactly summarised by introducing the 
functions
\be 
(A) \qquad 
f( \mu, \eta, \lambda)=e^{\mu a } e^{\eta v } e^{\lambda \tau \hspace{6em}} 
\ee
\be
(B) \qquad 
f( \mu, \lambda, \eta)=e^{\mu a} e^{\lambda \tau} e^{ \eta v \hspace{6em}} .
\ee

The choice of $ f( \mu, \eta, \lambda)$ in both  cases was dictated 
by simplicity
 of calculations.
It is now obvious that any element X of quantum Lie algebra is 
uniquely determined by the numerical function $f_x ( \mu, \eta, \lambda)$
 defined as
\be
f_x ( \mu, \eta, \lambda) \equiv <X,f( \mu, \eta, \lambda)>.
\ee
By appling the duality rules (4.1-4.4) and by multiple
use of Hausdorff formula
and some other tricks (cf Appendix) we arrive at the following quantum 
Lie algebra structures.

-The case (A)$\vspace{-1ex}$ 
\[
\vspace{-1ex} \Delta (H)=H \otimes I + I\otimes H \hspace{16em}
\]
\be \vspace{-1ex}
S(H)= - H \hspace{18em}
\ee
\[ \vspace{-1ex}
 \varepsilon (H)=0, \hspace{22em}
\]
\[
\Delta (K)=K \otimes I + e^{{-1\over \kappa}H} \otimes K -
{1 \over \rho} H e^{{-1\over \kappa}H} \otimes P \hspace{7em}
\]
\be
S(K)=-K e^{{-1\over \kappa}H}-{1\over\rho}HP e^{{-1\over \kappa}H}
 \hspace{10em} \ee
\[
 \varepsilon (K)=0, \hspace{22em}
\]
\[
 \Delta (P)=P \otimes I + e^{{-1\over \kappa}H} \otimes P \hspace{14em}
\]
\be
S(P)=-P e^{{-1\over \kappa}H} \hspace{16em}
\ee
\[
 \varepsilon (P)=0, \hspace{22em}
\]
\[
[H,P]= 0   \hspace{10em}
\]
\be
[K,P]= { -i\over{2\kappa}}P^2  \hspace{8em} 
\ee
\[
[K,H]=iP \hspace{9em}
\]
-the case (B) 
\[
\ba{c}
\Delta (H)=I \otimes H + H \otimes e^{-P \over \sigma}{\cosh(
{P\over \sqrt{\alpha \lambda}} )}  
-{ \alpha\over\sqrt{\alpha \lambda}}K\otimes e^{-P \over \sigma}
\sinh ({P \over \sqrt{\alpha \lambda}})   
\ea
\]
\be
\ba{c}
S(H)=-He^{P\over\sigma}\cosh({P\over\sqrt{\alpha \lambda}})-
{\lambda\over\sqrt{\alpha \lambda}}
Ke^{P\over\sigma}\sinh({P\over\sqrt{\alpha \lambda}}) \hspace{4em}
\ea
\ee
\[
 \varepsilon (H)=0, \hspace{23em} 
\]
\[
\ba{c}
\Delta(K)=I \otimes K + K \otimes e^{-P \over \sigma}{\cosh(
{P\over \sqrt{\alpha \lambda}} )}  
-{ \lambda\over\sqrt{\alpha \lambda}}H\otimes e^{-P \over \sigma}
\sinh ({P \over \sqrt{\alpha \lambda}}) 
\ea
\]
\be
\ba{c}
S(K)=-Ke^{P\over\sigma} \cosh({P\over\sqrt{\alpha \lambda}})-
{\alpha\over\sqrt{\alpha \lambda}}
He^{P\over\sigma}\sinh({P\over\sqrt{\alpha \lambda}}) \hspace{4em}
\ea
\ee
\[
 \varepsilon (K)=0, \hspace{23em}
\]
\[
\ba{c}
\Delta(P)=I \otimes P + P \otimes I \hspace{5em}
\ea
\]
\be
S(P)=-P \hspace{9em}
\ee
\[
 \varepsilon (P)=0,  \hspace{10em}
\]
\newpage
\[
[H,P]=0 \hspace{9em}
\]
\be
[K,P]=0 \hspace{9em}
\ee
\[
[K,H]={{i\sigma}\over -2}\left(e^{-2P\over \sigma } -1\right)
\hspace{3em} \]
and, in both cases, H,P, and K are hermitian. Some examples of actual
calculations are given in Appendix.
\zero
\section{The Lyakhovsky-Mudrov formalism}

	In order to find the quantum Lie algebras dual to our groups 
we can also use the formalism developed by Lyakhovsky and Mudrov [5], [6].
It is based on following theorem ( Lyakhovsky-Mudrov ):\\
\begin{em} 

Let $ \left\{ I, H_{1} ,...,H_{n} ,X_{1} ,...,X_{m} \right\} $ be
a basis of an associative algebra E over C verifying the conditions
\be
 \hspace{1em} \left[ H_{i} ,H_{j} \right] = 0, \hspace{10em}    i,j = 1,...,n.
\end{equation}
Let $\mu_{i} ,\nu_{j} (i,j=1,...,n) $ be a set of  $m \times n$ complex 
matrices such that
\be
[ \mu_{i} , \nu_{j} ] = [ \mu_{i} , \mu_{j} ] = [ \nu_{i} , \nu_{j} ] = 0,
   \hspace{6em}       i,j = 1,...,n. 
\ee

Let $\vec{X}$ be a vector with components $X_{l} (l = 1,...,m)$. \\
The coproduct  
\[ 
\Delta (I)=I \otimes I, \hspace{4em} 
 \Delta (H_i)=I \otimes H_i + H_i \otimes I
\]
\be
\Delta (\vec{X})= exp ( \sum_{i=1}^{n} {\mu_i}H_i ) \dot{ \otimes}
\vec{X} + \sigma \left( exp ( \sum_{i=1}^{n} {\nu_i}H_i ) 
 \dot{ \otimes} \vec{X} \right)
\ee
and the counit
\[
\varepsilon (I)=I , \hspace{2em}  \varepsilon (H_i)=0   \hspace{6em}
 i=1,...,n;  
\]
\be
 \varepsilon (X_l)=0 \hspace{12em} l=1,...,m; 
\ee

endow $ (E, \Delta, \varepsilon) $ with a coalgebra structure.
\end{em}

	With the help of this theorem we can find coalgebra structure.
 To this end we recall that the cocomutator $ \delta $ corresponds to the 
leading part of the co-antisymetic part of the coproducts
\be 
 \delta( \vec{X})= \Delta_{(1)} ( \vec{X})- \sigma \circ \Delta_{(1)}
( \vec{X}), \hspace{4em}
\ee
where
\be
\Delta_{(1)} (\vec{X})= ( \sum_{i=1}^{n} {\mu_i}H_i ) \dot{ \otimes}
\vec{X} + \sigma \left( ( \sum_{i=1}^{n} {\nu_i}H_i )
 \dot{ \otimes} \vec{X} \right)
\ee
 
is the first order term in all the parameters of (5.6).

Therefore matrices $ \mu_{i} $ and $ \nu_{i} $ can be determind from the
known form of $ \delta (X). $

	Now, if one is able to find a multiplication rules compatible with the 
coproduct one obtains a quantum algebra.
By applying this formalism to our case (which, actually, has been done in 
Ref.[6]) we arrive at the same form of coproduct as given by duality rules,
eq. (4.12 - 4.14; 4.16 - 4.18). 

Therefore, we have an alternative way to construct our quantum 
Galilei algebras.
\zero

\section{Acknowledgment}

The author acknowledges 
Prof. P. Ma\'slanka and
Prof. P. Kosi\'nski for a careful reading
of the manuscript and many helpful suggestions.
  Special thanks are also due to 
Prof. S. Giller, Dr. C. Gonera and Mrs A. Opanowicz
for valuable discussion.
\zero

\appendix
\section{Appendix}

We present some sample calculations concerning the dual
structures. Let us consider the (B) case,
\[
[a,v]={{i\over \alpha} \tau}-{i\over \sigma}v \hspace{6em} 
\]
\be
[a, \tau]={{i\over \sigma} \tau}+{i\over \lambda}v \hspace{6em} 
\ee
\[
[v, \tau]=0 \hspace{9em}
\]

$ \qquad \varepsilon^2 \sigma^2 = \alpha \lambda.  $

Let us calculate
\beqn
ff'=e^{\mu a} e^{\lambda \tau} e^{\eta v}
e^{\mu' a} e^{\lambda' \tau} e^{\eta' v} \hspace{8em} \nonumber \\ 
=e^{\mu a + \mu' a}{e^{\lambda{(e^{- \mu' a} \tau e^{\mu' a})}}}
{e^{\eta{(e^{- \mu' a} v e^{\mu' a})}}}{ e^{\lambda' \tau}}{e^{\eta'v}};
\eeqn
here prime means the second factor of tensor product
and the tensor product symbol $ \otimes $ has been omitted.
Denoting
\be
x( \mu')=e^{- \mu' a} \tau e^{\mu' a} \hspace{4em} x(0)= \tau
\ee
\be
y( \mu')=e^{- \mu' a} v e^{\mu' a} \hspace{4em} y(0)=v
\ee

we obtain the following differential equations

\be
\dot{x}(\mu')=e^{- \mu' a}[ \tau,a] e^{\mu' a} = {i\over \sigma}x(\mu')-
{i\over \lambda}y(\mu')
\ee
\be
\dot{y}(\mu')=e^{- \mu' a}[v,a] e^{\mu' a} = {i\over \sigma}y(\mu')-
{i\over \alpha}x(\mu')
\ee
or, in matrix form
\be 
\stackrel{\bullet}{ \left(
\begin{array}{c}
x(\mu') \vspace{1ex} \\y(\mu')
\end{array}
\right)}=
 \left(
\begin{array}{cc}
{i\over \sigma} & {-i\over \lambda} \vspace{1ex} \\
{-i\over \alpha} &{i\over \sigma}
\end{array}
\right)
 \left(
\begin{array}{c}
x(\mu') \vspace{1ex} \\y(\mu')
\end{array}
\right).
\ee

Thus the  solution to eq.(A.7) reads

\be 
 \left(
\begin{array}{c}
x(\mu') \vspace{1ex}\\y(\mu')
\end{array}
\right)=
e^{i \mu' A}
 \left(
\begin{array}{c}
\tau \vspace{1ex}\\ v
 \end{array} 
 \right) \hspace{5em}
\ee
where 
\[
A=
\left(
\begin{array}{cc}
{1\over \sigma} & {-1\over \lambda} \vspace{1ex} \\
{-1\over \alpha} & {1\over \sigma}
\end{array} 
\right). \hspace{8em}
\]
It is easy to check that
\be
e^{i \mu' A}=e^{i \mu' \over \sigma}
\left(
\begin{array}{cc}
{\cosh\left({ \mu'\over \sqrt{\alpha \lambda}} \right)} & 
{- \alpha\over\sqrt{\alpha \lambda}}
{\sinh \left({\mu' \over \sqrt{\alpha \lambda}} \right)} \vspace{3ex}\\
{- \lambda\over\sqrt{\alpha \lambda}}
{\sinh \left({\mu' \over \sqrt{\alpha \lambda}} \right)} &
{\cosh { \left(\mu'\over \sqrt{\alpha \lambda} \right)}} 
\end{array}
\right)
\ee 
and, consequently  
\be
\begin{array}{c}
x= e^{i \mu' \over \sigma}{ \left(\tau{\cosh(
{ \mu'\over \sqrt{\alpha \lambda}} )} - 
{ \alpha\over\sqrt{\alpha \lambda}}
{v \sinh ({\mu' \over \sqrt{\alpha \lambda}})} \right)} \vspace{2ex} \\
y= e^{i \mu' \over \sigma}{ \left({- \lambda\over\sqrt{\alpha \lambda}}
{\tau \sinh({ \mu' \over \sqrt{\alpha \lambda}})}+
{v \cosh ( {\mu'\over \sqrt{\alpha \lambda}})} \right).}
\end{array}
\ee
Therefore
\be
\begin{array}{c} 
ff'=\exp {( \mu + \mu')a}  \hspace{18em} \vspace{2ex} \\ 
\cdot \exp {\left[ e^{i \mu' \over \sigma}{ \left(
{- \lambda\over\sqrt{\alpha \lambda}}
{\eta \sinh({ \mu' \over \sqrt{\alpha \lambda}})}+
{\lambda \cosh ( {\mu'\over \sqrt{\alpha \lambda}})} \right)}
+ \lambda'\right] \tau} \vspace{2ex} \\
\cdot \exp {\left[e^{i \mu' \over \sigma}{ \left(\eta{\cosh(
{ \mu'\over \sqrt{\alpha \lambda}} )} - 
{ \alpha\over\sqrt{\alpha \lambda}}
{\lambda \sinh ({\mu \over \sqrt{\alpha \lambda}})} \right)} 
+ \eta'\right]v}
\end{array}
\ee  
and from this formula we obtain 
\be
\ba{c}
\Delta (H)=I \otimes H + H \otimes e^{-P \over \sigma}{\cosh(
{P\over \sqrt{\alpha \lambda}} )}  
-{ \lambda\over\sqrt{\alpha \lambda}}K\otimes e^{-P \over \sigma}
\sinh ({P \over \sqrt{\alpha \lambda}}) \vspace{2ex} \\ 
\Delta(K)=I \otimes K + K \otimes e^{-P \over \sigma}{\cosh(
{P\over \sqrt{\alpha \lambda}} )}  
-{ \alpha\over\sqrt{\alpha \lambda}}H\otimes e^{-P \over \sigma}
\sinh ({P \over \sqrt{\alpha \lambda}}) \vspace{2ex} \\
\Delta(P)=I \otimes P + P \otimes I. \hspace{18em}
\ea
\ee
On the other hand
\be
\Delta (f) =e^{\mu (a+a' +v \tau')} e^{ \lambda (\tau + \tau')}
 e^{ \eta (v+v')}.
\ee
In order to calculate $ \Delta (f) $ we use the following trick:
we write 
\be
\Delta (f)= e^{ \mu (a+a')} \chi e^{ \lambda( \tau + \tau')}
e^{ \eta(v+v')}
\ee where \be
\chi=e^{- \mu(a+a')} e^{ \mu(a+a'+v \tau')} \hspace{3em}
\chi(0)=1. \ee
Again differentiating both sides
with respect to $ \mu $ we got
\begin{equation}
 \dot{ \chi} =( e^{-{ \mu}a} v e^{ \mu a} ) 
 ( e^{-{ \mu}a'} \tau' e^{ \mu a'} ) \chi \hspace{5em}
\end{equation}
and
\be
\ba{c}
{{\dot{ \chi}}\over{ \chi}} = e^{2i \mu \over \sigma}
{ \left(v{\cosh({ \mu\over \sqrt{\alpha \lambda}})} - 
{ \lambda\over\sqrt{\alpha \lambda}}
{\tau \sinh ({\mu \over \sqrt{\alpha \lambda}})} \right)} \vspace{3ex} \\
\hspace{3em} {\left( \tau'{\cosh ( {\mu\over \sqrt{\alpha \lambda}})} 
-{ \alpha\over\sqrt{\alpha \lambda}}
{v' \sinh({ \mu \over \sqrt{\alpha \lambda}})}\right).}
\end{array}
\ee
The terms $ vv' $ and $ \tau \tau' $ do not contribute to the product
of different generators, so they can be neglected in what follows.
Therfore, up to irrelevant terms 
\be 
\begin{array}{l}
\Delta ( f)=e^{\mu a} e^{\mu a'} e^{\lambda \tau}
e^{ \lambda \tau'} e^{\eta v} e^{\eta v'} 
\hspace{19em} \vspace{1ex} \\
e^{ \tau v' \left[ {-i \sqrt{\alpha \lambda}}
\left( 
{{e^{({{2i \mu} \over \sigma})}} \over{2{(1-{{\alpha \lambda}\over \sigma^2})}}}
\left(
\sinh({{i \mu} \over \sqrt{\alpha \lambda}})
\left(
\cosh({{i \mu} \over \sqrt{\alpha \lambda}}) 
+{ \sqrt{\alpha \lambda}\over \sigma}
\sinh({{i \mu} \over \sqrt{\alpha \lambda}})
\right)
-{\sigma \over{2\sqrt{\alpha \lambda}}}
\right)
\right) +
{ \sigma \over{ 4 \sqrt{\alpha \lambda}
 {(1- {{ \alpha \lambda} \over \sigma^2})}}} \right] } \vspace{1ex}
 \nonumber  \\
e^{ \tau' v \left[ {-i \sqrt{\alpha \lambda}}
\left( 
{{e^{({{2i \mu} \over \sigma})}} \over{2{(1-{{\alpha \lambda}\over \sigma^2})}}}
\left(
\cosh({{i \mu} \over \sqrt{\alpha \lambda}})
\left(
\sinh({{i \mu} \over \sqrt{\alpha \lambda}}) 
+{ \sqrt{\alpha \lambda}\over \sigma}
\cosh({{i \mu} \over \sqrt{\alpha \lambda}})
\right)
+{\sigma \over{2\sqrt{\alpha \lambda}}}
\right)
\right) -
{{ \sigma(1-2{{\sqrt{\alpha \lambda}\over \sigma^2 }})}
\over{ 4 \sqrt{\alpha \lambda}
{(1- {{ \alpha \lambda} \over \sigma^2})}}} \right] } \nonumber  
\ea
\ee. 
Hence, the duality rule (4.1) gives 
\[
[H,P]=0 \hspace{6em}
\]
\be
[K,P]=0 \hspace{6em}
\ee
\[
[K,H]={{i\sigma}\over 2}\left(1- e^{-2P\over \sigma }\right).
\]
The antipode and counit can be obtained either from the relevant
duality relations or directly from the Hopf algebra rules:
\be
(id \otimes \varepsilon) \Delta = id
\ee
\be
m \circ (id \otimes S) \Delta = \varepsilon.
\ee
In this way we have
\be
\varepsilon (X)=0 \hspace{17em} \vspace{2ex}
\ee
\be
\ba{c}
S(H)=-He^{P\over\sigma}\cosh({P\over\sqrt{\alpha \lambda}})-
{\lambda\over\sqrt{\alpha \lambda}}
Ke^{P\over\sigma}\sinh({P\over\sqrt{\alpha \lambda}}) \vspace{2ex}
\\
S(K)=-Ke^{P\over\sigma}\cosh({P\over\sqrt{\alpha \lambda}})-
{\alpha\over\sqrt{\alpha \lambda}}
He^{P\over\sigma}\sinh({P\over\sqrt{\alpha \lambda}}) \vspace{2ex}\\
 S(P)=-P; \hspace{16em}
\ea
\ee
moreover coproduct is homomorphism of this algebra i.e.
\[
\Delta([ , ])= [ \Delta, \Delta ].
\]
Let us now check eq(4.4):
\beqn
S^{-1}\left(\left(e^{\lambda \tau} 
e^{\mu a} e^{\eta v}\right)^*\right)= 
S^{-1}\left(e^{\eta^* v} 
e^{\mu^* a} e^{\lambda^* \tau}\right)= \vspace{2ex} \nonumber\\
e^{-\lambda^* \tau} 
e^{\mu^*(-a +\tau v )} e^{-\eta^* v}. \hspace{10em} 
\eeqn 
\beqn
\left[ S\left(\left(e^{\lambda \tau} 
e^{\mu a} e^{\eta v}\right)^*\right)\right]^*= 
\left[ S\left(e^{\eta^* v} 
e^{\mu^* a} e^{\lambda^* \tau}\right) \right]^*=
 \vspace{1ex}\hspace{3em} \nonumber\\
\left(e^{-\lambda^* \tau} 
e^{\mu^*(-a +\tau v )} e^{-\eta^* v}\right)^*= 
e^{-\eta v} e^{\mu(-a+v\tau)} e^{-\lambda \tau}= 
\hspace{2em} \nonumber\\
 S^{-1}\left(e^{\lambda \tau} 
e^{\mu a} e^{\eta v}\right). \hspace{15em}
 \eeqn 
By applying the formula (4.3) we get
\beqn
<H^*,e^{\lambda \tau} e^{\mu a} e^{\eta v}>=
<H,e^{-\lambda^*\tau} e^{\mu^*(-a+\tau v)} e^{-\eta^*v}>^* =
\\ 
(-i\lambda^*)^*=i\lambda
=<H,e^{\lambda \tau} e^{\mu a} e^{\eta v}> \hspace{8em} \nonumber
\eeqn
hence
\[
H^*=H; \hspace{12em}
\]
In this same way one obtains \(
K^*=K, \hspace{1em}
\)
\(
P^*=P. \hspace{1em}
\)
\zero



\begin{thebibliography}{99}
\newcommand{\byauthors}[1]{#1 }
\newcommand{\journal}[1]{ #1 }
\newcommand{\reftitle}[1]{{\it #1} }
\newcommand{\volumin}[1]{{\bf #1} }
\newcommand{\eref}{.}
\newcommand{\refyear}[1]{(#1) }
\bibitem{1}\byauthors{E. Kowalczyk }
\journal{Acta Physica Polonica B}
\volumin{28}
\refyear{1997}\eref
\bibitem{2}\byauthors{V.Chan, A. Pressley }
\journal{Cambrdige Univ.Press}
\refyear{1994}
\eref
\bibitem{3}\byauthors{S. Zakrzewski }
\journal{Comm. Math. Phys.}
\volumin{185}
\refyear{1997}
\eref
\bibitem{4}\byauthors{P.Ma\'slanka }
\journal{J.Math.Phys}
\volumin{35}
\refyear{1994}
\eref
\bibitem{5}\byauthors {A. Ballesteros, F. J. Herranz, 
P. Parashar}
\journal{J. Phys. A: Math. Gen.}
\volumin{30}
\refyear{1997}
\eref
\bibitem{6}\byauthors {A. Ballesteros, F. J. Herranz} 
\reftitle{Lie Bialgebra Quantizations of the Oscillator 
Algebra and their Universal {\em R}-Matrices }
\journal{XXI ICGTMP, Goslar (Germany)}
\refyear{1996}
\eref
\end{thebibliography}
\end{document}